\documentclass[12pt,leqno,twoside]{amsart}
\usepackage{amssymb,amsfonts,amsmath,amsthm}
\usepackage{latexsym}
\usepackage{verbatim}
\usepackage{epsfig}

\newtheorem{theorem}{Theorem}[section]
\newtheorem{corollary}[theorem]{Corollary}

\newtheorem{proposition}[theorem]{Proposition}
\theoremstyle{remark}
\newtheorem{remark}[theorem]{\sc Remark}
\theoremstyle{remark}

\theoremstyle{definition}
\newtheorem{definition}[theorem]{Definition}
\theoremstyle{remark}
\newtheorem{example}[theorem]{\sc Example}

\theoremstyle{remark}

\theoremstyle{remark}

\numberwithin{equation}{section}

\renewcommand{\Box}{_\square}    %\diamond

\newcommand{\cal}{\mathcal}

\newcommand{\diff}{{\rm{diff}}}
\renewcommand{\int}{{\rm{int}}}
\newcommand{\gen}{{\rm{gen}}}

\newcommand{\Sing}{{\rm{Sing\hspace{2pt}}}}
\newcommand{\Crit}{{\rm{Crit}}}
\newcommand{\Disc}{{\rm{Disc\hspace{2pt}}}}

\newcommand{\ity}{{\infty}}

\newcommand{\e}{\varepsilon}
\newcommand{\fin}{\hspace*{\fill}$\Box$\vspace*{2mm}}
\newcommand{\fisi}{{$\mathcal F${\rm{ISI}\hspace{3pt}}}}

\newcommand{\cB}{{\cal B}}

\newcommand{\cW}{{\cal W}}

\newcommand{\bC}{{\mathbb C}}

\newcommand{\bP}{{\mathbb P}}

\newcommand{\bX}{{\mathbb X}}
\newcommand{\bY}{{\mathbb Y}}

\begin{document}

\title[Singularity exchange]
 {Singularity exchange at the frontier of the space}

\author{\sc Dirk Siersma}

\address{Mathematisch Instituut, Universiteit Utrecht, PO
Box 80010, \ 3508 TA Utrecht
 The Netherlands.}

\email{siersma@math.uu.nl}

\author{\sc Mihai Tib\u ar}

\address{Math\' ematiques, UMR 8524 CNRS,
Universit\'e de Lille 1, \  59655 Villeneuve d'Ascq, France.}

\email{tibar@math.univ-lille1.fr}

%\thanks{}



%\date{\today}

%\dedicatory{}

%\commby{}

%%% ----------------------------------------------------------------------

\begin{abstract}
In deformations of polynomial functions one may encounter
 ``singularity exchange at
infinity'' when singular points  disappear from the space and produce
  ``virtual'' singularities which have an influence on the topology
 of the limit polynomial.  We find several rules 
of this exchange phenomenon, in which the total quantity of singularity
 turns out to be not conserved in general. 

\end{abstract}

%%% ---------------------------------------------------
\maketitle
%%% ---------------------------------------------------

\setcounter{section}{0}
\section{Introduction}

More than 20 years ago, Broughton \cite{Br} pioneered the study of  the topology of
the fibres of complex polynomial functions of a certain type (called B-type
here). Even if the study of the singular fibration produced by a polynomial in
the affine space took a certain ampleness ever since, 
 {\em families} of such polynomial
functions have been
considered only sporadically. We study here families of
 polynomial functions by focussing on the transformation of
singularities in the neighbourhood of infinity, a phenomenon which we already
 remarked in \cite{ST-gendefo}.
 This is a natural and chalenging topic inside mathematics since the atypical
fibres of a polynomial turn out to be not only due to the ``visible''
singularities, but also to the ``bad'' asymptotic behaviour at the infinite
frontier of the space.
 We deal here with the evolution and interaction of
singularities in deformations at the infinite frontier of the space, in what concerns
the phenomenon of conservation or
non-conservation of certain  numbers attached to singularities (that we recall
below).

Let $\{ f_s\}_s$ be a holomorphic family of complex polynomial functions $f_s
: \bC^n \to \bC$, for $s$ in a small neighbourhood of $0\in \bC$.
  For a fixed polynomial function $f_s$ there is a well
 defined general fibre $G_s$, since the set of atypical values $\Lambda(f_s)$
 is a finite set. When specialising to $f_0$,
  the number of atypical
 values may vary (decrease, increase or be constant) and the topology of
 the general fibre may change.
 We consider
 constant degree families within certain classes
of polynomials  (F-class $\subset$ B-class
$\subset$ W-class, cf  Definition \ref{d:broughton})
which have the property that the {\em vanishing cycles of $f_s$}
 (i.e. the generators of the reduced homology of $G_s$) are
 concentrated in dimension $n-1$ and are localisable at finitely many points,
 in the affine space or in the part at infinity of the projective
 compactification of some fibre of $f_s$.
In the affine space $\bC^n$, such a point is a singular
point of $f_s$. The sum of all affine Milnor numbers
 is the total Milnor number $\mu(s)$, which has an algebraic interpretation as
 the dimension of the quotient algebra $\bC[x_1, \ldots , x_n]/(\frac{\partial
 f}{\partial x_1}, \ldots ,\frac{\partial f}{\partial x_n})$.
  Singularities at infinity are equipped with so-called Milnor-L\^e
 numbers (cf \cite{ST}) and their sum is denoted by $\lambda(s)$. Then the
 Euler characteristic of  the generic fiber $G_s$ is $1+
 (-1)^{n-1}(\mu(s)+\lambda(s))$.

  A natural problem which arises is
 to understand the behaviour, when $s\to 0$, of the $\mu$ and $\lambda$-singularities, 
which support the vanishing
cycles of $f_s$.
 It is well known and easy to see that, for singularities which tend to 
a $\mu$-singular point,
  the total number of local vanishing cycles is constant, in other words the
local balance law is {\em conservative}. However some
 $\mu$-singularities may tend to infinity  and change into
 $\lambda$-singularities; this is the phenomenon we address here. 
First, in full generality, for any deformation, we get the:
\begin{itemize}
  \item[$-$] {\em global lower semi-continuity} of the highest Betti number:
 $b_{n-1}(G_0) \le b_{n-1}(G_s)$ (Proposition \ref{t:spec}).
\end{itemize}

Next, focussing on constant degree deformations inside the B-class, we prove
 several facts on the singularity exchange at
infinity:

\begin{itemize}
  \item[$-$] the number of local vanishing cycles
  of $\mu$ and $\lambda$-singularities tending to a $\lambda$-singular point
 is {\em lower semi-continuous}, but it is not conserved in general 
(Theorem \ref{t:count}).

 \item[$-$]
 in $(\mu + \lambda)$-constant deformations, the local balance law at any
$\lambda$-singularity of $f_0$ is {\em conservative} and {\em atypical values cannot escape to infinity} (Corollary \ref{c:noloss}).

\item[$-$]
 in $(\mu + \lambda)$-constant deformations,
the monodromy fibrations over any admissible loop (in particular,
 the monodromy fibrations at infinity) are isotopic in the family,
 whenever $n\not=3$ (Theorem \ref{t:mono}).

\item[$-$]
 in deformations  with constant generic singularity type at infinity,
 $\lambda$-singularities of $f_0$ are {\em locally persistent} in $f_s$
but cannot split such that more than one
 $\lambda$-singularity occurs in the same fibre (Theorem \ref{t:inf}).

\item[$-$]
 in deformations inside the F-class,  a $\lambda$-singularity cannot be deformed into only
 $\mu$-singularities (Corollary \ref{c:lambda}).

\end{itemize}

The semi-continuity results (first two of the above list) are certainly
related to the semi-continuity of the {\em spectrum}, 
a result proved by N\'emethi and Sabbah \cite{NS} for the class
of ``weakly tame'' polynomials. Their class excludes by definition the
$\lambda$-singularities, but on the other hand the spectrum
 (defined with Hodge theoretical ingredients)
gives more refined information than the total  Milnor number.
 It is also interesting to remark that the lower semi-continuity in all these 
results is opposite to the
{\em upper semi-continuity}  in  case
of deformations of  holomorphic function germs. 

We end by supplying with a zoo of examples which illustrate various aspects of 
the exchange phenomenon.

%%%%%%%%%%%%%%%

%%%%%%%%%%%%%%%
%%%%%%%%%%%%%%%%%%%%%%%%%
\section{Deformations in general}\label{ss:defo}

%\marginpar{fffIt is well known}

It is well known that
  the $(n-1)$th Betti number of the Milnor fibre of a holomorphic function
  germ
 is upper semi-continuous, i.e. it does not decrease under
specialisation.
In case of a polynomial $f_s : \bC^n \to \bC$, the role of the Milnor
fibre is
played by the general fibre $G_s$ of $f_s$. This is a Stein
manifold of
dimension $n-1$ and therefore it has the homotopy type of a CW complex of
 dimension $\le n-1$, which is also finite,
since $G_s$ is algebraic. Moreover, the $(n-1)$th homology group with integer
 coefficients is free.
We prove the following general specialisation result.

%%%%%%%%%%%%%%%%%%%%%%%%%%
%%%%%%%%%%%%%%%%%%%%%%%%%%%%%%%%%%%%%%%%%%%%%%%%%%
\begin{proposition}\label{t:spec}
Let $P : \bC^n \times \bC^k \to \bC$ be any holomorphic deformation of a
 polynomial $f_0:= P(\cdot, 0) : \bC^n
\to \bC$. Then
 the general fibre $G_0$ of $f_0$ can be naturally embedded into the
general
fibre $G_s$ of $f_s$, for $s\not= 0$ close enough to $0$.
The embedding $G_0\subset G_s$ induces an inclusion $H_{n-1}(G_0)
 \hookrightarrow  H_{n-1}(G_s)$ which is compatible with the intersection form.
\end{proposition}

\begin{proof}
It is enough to consider a
$1$-parameter family of hypersurfaces $\{ f_s^{-1}(t)\}_{s\in L}\subset
\bC^n$, for fixed $t$,
where $L$ denotes some parametrised complex curve through $0$. We denote by
$X_t$ the total
space over a small neighbourhood $L_\varepsilon$ of $0$ in $L$.
By choosing $t$ generic enough,  we may assume that $f_s^{-1}(t)$ is a
generic
fibre of $f_s$, for $s$ in a small enough neighbourhood of $0$. Let $\sigma :
X_t \to L_\varepsilon$ denote the projection. Now $X_t$ is the total space of a
 family of non-singular hypersurfaces.
 Since $\sigma^{-1}(0)$ is an affine
hypersurface, by taking a large enough radius $R$, we get $\partial \bar
 B_{R'}
 \pitchfork
\sigma^{-1}(0)$, for all $R'\ge R$. Moreover, the sphere $\partial \bar B_R$
 is
 transversal to all
nearby fibres $\sigma^{-1}(s)$, for small enough $s$.
 It follows that the projection $\sigma$ from the pair of spaces $(X_t\cap
(B_R\times \bC), X_t\cap
(\partial \bar B_R\times \bC))$ to $L_\varepsilon$
 is a proper submersion
and hence, by Ehresmann's
theorem, it is a trivial fibration. By the above transversality argument,
we have  $B_R \cap\sigma^{-1}(0)\stackrel{\diff}{\simeq} B_R
\cap\sigma^{-1}(s)$. This shows the first claim.

  The affine hypersurfaces $\sigma^{-1}(s)$ are finite cell complexes
of dimension $\le n-1$. By the classical Andreotti-Frankel \cite{AF} argument
 for the distance
function, the hypersurface
 $\sigma^{-1}(s)$ is obtained from $B_R \cap\sigma^{-1}(s)$ by adding cells
of index at most $n-1$. This shows that $H_n(G_s, G_0) = 0$, so the second
claim. The compatibility with the intersection form is standard.
\end{proof}
%%%%%%%%%%%%%%
Under certain conditions we can also compare
the ``monodromy fibrations at infinity'' in the family, see \S\ref{ss:mon}.
Proposition \ref{t:spec}  will actually be exploited through the
{\em semi-continuity of the highest
 Betti number}, as a consequence of the inclusion of homology groups:
%%%%%%%%%%%%%%%%%%%%%%%%%%
\begin{equation}\label{c:globsemi}
b_{n-1}(G_s) \ge b_{n-1}(G_0), \mbox{ for } s \mbox{ close enough to } 0.
\end{equation}
%%%%%%%%%%%%%%%%%%%

 %%%%%%%%%%%%%%%%%%%%%%%%%%%%%%%%%%%%%
\section{Compactification of families of polynomials}\label{s:class}

%%%%%%%%%%%%%%%

We shall now focus on polynomials
for which the singularities at infinity are isolated, in a sense that we make precise here.

Let $P$ be a deformation of $f_0$, i.e. $P: \bC^n \times \bC^k \to \bC$  is a
 family of polynomial functions
$P(x,s) = f_s(x)$ such that $f_0 = f$.
We assume in the following that
our deformation depends holomorphically on the parameter $s\in \bC^k$. We also
 assume
 that $\deg f_s$ is independent on $s$, for $s$ in some neighbourhood of $0$,
 and
 we denote it by $d$.
We
attach to $P$ the following hypersurface:
\[ \bY = \{ ([x:x_0], s, t) \in \bP^n \times \bC^k \times \bC \mid \tilde
P(x,x_0,s) - tx_0^d = 0 \},\]
where $\tilde P$ denotes the homogenized of $P$ by the variable $x_0$,
considering $s$ as parameter varying in a small neighbourhood of $0\in \bC^k$.
Let $\tau : \bY \to \bC$ be the projection to the $t$-coordinate. This
extends the map $P$ to a proper one in the sense that $\bC^n\times \bC^k$ is
 embedded in $\bY$ (via the graph of $P$) and 
 $\tau_{|\bC^n\times \bC^k} = P$. Let $\sigma : \bY \to
\bC^k$ denote the projection to the $s$-coordinates.

\noindent {\bf Notations.}
$\bY_{s,*} := \bY\cap \sigma^{-1}(s)$, $\bY_{*, t} := \bY\cap
\tau^{-1}(t)$ and
$\bY_{s,t} := \bY_{s,*}\cap \tau^{-1}(t) = \bY_{*, t}\cap
\sigma^{-1}(s)$. Note
that $\bY_{s,t}$ is
the closure in $\bP^n$ of the affine hypersurface $f_s^{-1}(t)\subset \bC^n$.

Let $\bY^\ity := \bY\cap \{ x_0=0\} = \{
P_d(x,s) = 0\}\times \bC$ be the hyperplane at infinity of $\bY$, where $P_d$
is the degree $d$ homogeneous part of $P$ in variables $x\in \bC^n$.
Remark that for any fixed $s$, $\bY^\ity_{s,t} :=\bY_{s, t}\cap \bY^\ity$
does not depend on $t$.

%%%%%%%%%%%%%%%%%
%%%%%%%%%%%%%%%%%%%%%%%%%%%%%%%%%%%%%
\begin{definition}\label{d:broughton}
We consider the following classes of polynomials:
\begin{itemize}
\item[(i)] $f$  is a  {\em F-type} polynomial if its compactified fibres
and their restrictions to the hyperplane at infinity
have at most isolated singularities.
\item[(ii)] $f$  is a {\em B-type} polynomial if its compactified fibres
have at most isolated singularities.
 \end{itemize}
\end{definition}
%%%%%%%%%%%%%%%%
It follows that F-class $\subset$ B-class. They are both contained into the W-class, which consists
polynomials for which the proper  extension $\tau : \bX \to \bC$ has
  only isolated singularities with respect to some Whitney stratification of
  $\bX$ such that $\bX^\ity$ is a union of strata, see \cite{ST}. The notation
  $\bX$ stands for $\bY$ when a single polynomial is considered (i.e. there is
  no parameter $s$).

 In two variables, if $f$ has
isolated singularities in $\bC^2$, then it is automatically of F-type.  Deformations inside the
F-class were introduced in
 \cite{ST-gendefo} under the name \fisi deformations.
   Broughton \cite{Br}
 considered for the first time B-type polynomials  and studied the
topology of their general fibers.  The W-class
of polynomials appears in \cite{ST}. In deformations of a polynomial
 $f_0$ we usually require
to stay inside the same class but we may also deform into a "less singular" class
(like B-type into F-type, Example \ref{e:1}).

%%%%%%%%%%%%%%%%%

The singular locus of $\bY$,
$ \Sing \bY := \{ x_0 =0, \frac{\partial P_d}{\partial x}(x,s) =0, \
 P_{d-1}(x,s)=0,
\frac{\partial P_d}{\partial s}(x,s)=0 \} \times \bC$
is included in $\bY^\ity$ and is a product-space by the $t$-coordinate.
It depends only on the degrees $d$ and $d-1$ parts of $P$ with respect to the
 variables $x$.

Let $\Sigma := \{ x_0 =0, \frac{\partial P_d}{\partial x}(x,s) =0, \
 P_{d-1}(x,s)=0\}\subset \bP^{n-1}\times \bC^k$. If we fix $s$, the singular locus of
$\bY_{s,*}$ is the  analytic set $\Sigma_s \times \bC$, where
$\Sigma_s := \Sigma \cap \{ \sigma =s\}$, and it is the union of the
 singularities at the hyperplane at
infinity of
the hypersurfaces $\bY_{s,t}$, for $t\in \bC$.

We denote by  $W_s := \{ [x] \in \bP^{n-1} \mid
\frac{\partial P_d}{\partial x}(x,s) =0 \}$ the set of points at infinity
where
 $\bY^\ity_{s,t}$ is singular, in other words where $\bY_{s,t}$ is either
 singular or tangent to $\{ x_0 =0\}$. It does not depend on $t$ and we
 have $\Sigma_s\subset W_s$.

%%%%%%%%%%%%%%%%%%%%%%%%%%%%%%%%%%%%%
\begin{remark}\label{r:broughton}
From the above definition and the expressions of the singular loci we have the
 following characterisation:
\begin{itemize}
\item[(i)]  $f_0$  is a B-type polynomial $\Leftrightarrow$ $\dim \Sing f_0 \le 0$
and $\dim \Sigma_0 \le 0$,
\item[(ii)]  $f_0$  is a F-type polynomial $\Leftrightarrow$  $\dim \Sing f_0 \le 0$
and  $\dim W_0 \le 0$.
\end{itemize}

Let us also remark that $\dim \Sigma_0 \le 0$ (respectively $\dim W_0 \le 0$)
 implies
 that $\dim \Sigma_s \le 0$ (respectively $\dim W_s \le 0$), whereas
 $\dim
 \Sing f_0 \le 0$ does not imply automatically $\dim \Sing f_s \le 0$ for
 $s\not=0$.
\end{remark}

%%%%%%%%%%%%%%%%%%%%%
%%%%%%%%%%%%%%%%%%%%%%

%%%%%%%%%%%%%%%%%%%%%

%\section{Local semi-continuity}\label{s:local}

%%%%%%%%%%%%%%%%%%%%%%%%%

\section{Semi-continuity at infinity}\label{ss:semi}

Let $P$ be a  deformation of $f_0$ such that $f_s$ is of
W-type, for all $s$ close enough to $0$.
It is shown in \cite{Pa, ST} that the vanishing cycles of $f_s$ (for fixed $s$)
 are concentrated in dimension $n-1$ and are localized at well-defined points, 
either in the affine space  or at infinity.
 We shall call them {\em $\mu$-singularities} and {\em $\lambda$-singularities}
 respectively. To such a singular point $p\in\bY_{s,*}$ one associates its
 local Milnor number denoted $\mu_{p}(s)$ or its Milnor-L\^e number
 $\lambda_{p}(s)$.
Let $\mu(s)$ be the {\em total Milnor
number}, respectively $\lambda(s)$ be the {\em total Milnor-L\^e number at
 infinity}, where  $b_{n-1}(G_s) = \mu(s) + \lambda(s)$.

By \cite{ST}, the atypical fibers of a W-type polynomial $f_s$ are exactly
 those fibers which contain $\mu$ or  $\lambda$-singularities; equivalently, those of
 which the Euler characteristic is different from
 $\chi(G_s)$. We denote by $\Lambda(f_s)$ the set
 of
 atypical values of $f_s$.

 The above cited facts
together with our semi-continuity result (\ref{c:globsemi})
show that, for $s$ close to $0$
 we have:
\[ \mu(s) + \lambda(s) \ge \mu(0) + \lambda(0).\]

%%%%%%%%%%%%%%
%%%%%%%%%%%%%%%%%%%%%
\begin{remark}\label{r:ab}
The total Milnor number
 $\mu(s)$ is lower
 semi-continuous under specialization $s\to 0$.
In case $\mu(s)$ decreases, we say that there is {\em loss of
 $\mu$ at infinity},
since this may only happen when one of the two following phenomena occur:
\begin{enumerate}
 \item\label{i} the modulus of some critical point tends to infinity
and the corresponding critical value is bounded (\cite[Example 8.1]{ST-gendefo});
\item\label{ii} the modulus of some critical value tends to infinity
 (\cite[Examples (8.2) and (8.3)]{ST-gendefo}).
\end{enumerate}
 In contrast to $\mu(s)$, it turns out that $\lambda(s)$ is not
 semi-continuous; under specialization,
 it can increase
  or decrease (Example \ref{e1:a}, \ref{e1:c}).
 Moreover, the $\lambda$-values may behave like the critical values in case (b) above, see
 Example \ref{e:1}.
\end{remark}
%%%%%%%%%%%%%%%%%%%%%%%%%%%%%%%%%%%%%%%%%%%%
 To understand the behaviour of $\lambda(s)$ in more detail,
 we focus on the B-class. The following result extends our
 \cite[Theorem 5.4]{ST-gendefo} and needs a more involved proof, which will be
 given in \S \ref{s:proof}.
%%%%%%%%%%%%%%%%%%%%%%%%%%%%%%%%%%%%%%%
\begin{theorem}\label{t:count}{\bf (Lower semi-continuity at
 $\lambda$-singularities)}\\
Let $P$ be a constant degree one-parameter deformation inside the B-class. Then,
locally at
any $\lambda$-singularity $p\in \bY_{0,t}$ of $f_0$, we have:
\[ \lambda_p(0) \le \sum_i \lambda_{p_i}(s) + \sum_j \mu_{p_j}(s),\]
where $p_i$ are the $\lambda$-singularities and $p_j$ are the $\mu$-singularities of $f_s$ which
tend to the point $p$ as $s\to 0$.
\end{theorem}

%%%%%%%%%%%%%%%%%%%%%%%

%%%%%%%%%%%%%%%%%%%%%%%%%%%%
%%%%%%%%%%%%%%%%%%
%%%%%%%%%%%%%%
%%%%%%%%%%%%%%%%%%

\section{Persistence of $\lambda$-singularities}\label{persistance}

%%%%%%%%%%%%%%%%%%%%%%%%%%

 In order to get further information on the $\mu\mapsto\lambda$ exchanges
 we focus on two sub-classes of the B-class.
 In this section we define cgst-type deformations and in the next section we
 study deformation with constant $\mu +\lambda$.

Let us first remark that for a deformation $\{ f_s\}_s$ inside the B-class the compactified fibres
of $f_s$ have only isolated singularities. The positions of these singularities depend only on $s$
(and not on $t$). When $s\to 0$ these singularities can split or disappear.

Let us take some $x(0)\in \Sigma_0$. Take $t\not\in \Lambda(f_0)$ and assume without dropping
generality that $t\not\in \Lambda(f_s)$ for all small enough $s$. Lazzeri's
non-splitting argument, see \cite{La} and also
\cite{AC,Le}, tells us that the Milnor number of
 $\bY_{0,t}$ at $(x(0),t)$ is strictly larger than  the
sum of the Milnor numbers of $\bY_{s,t}$ at all points $(x(s),t)\in \Sigma_s \times \{t\}$ such
that $x_i(s) \to x(0)$, unless there is only one such singular point $x(s)$ and the Milnor number
of $\bY_{s,t}$ at $(x(s),t)$ is independent on $s$. In the latter case we say that {\em the cgst
assumption holds}.

\begin{definition}\label{d:type}
We say that a constant degree deformation  inside the B-class has {\em
constant  generic singularity type at
 infinity at some point $x(0)\in \Sigma_0$} if the  cgst assumption holds.

If the cgst assumption holds at all points in $\Sigma_0$, then
we simply say {\em constant generic singularity type at infinity}.
\end{definition}

Note that in the B-class the cgst assumption does not imply that $b_{n-1}(G_s)$ is constant not
necessarily constant in the B-class (see Example \ref{e:1}). We send to Remark \ref{r:rel} for
further comments on cgst.

\begin{theorem}\label{t:inf}
Let $P$ be a constant degree deformation, inside the B-class,
with constant generic singularity type at
 infinity. Then:
\begin{enumerate}
\item  $\lambda$-singularities of $f_0$ are locally persistent in $f_s$.
\item a $\lambda$-singularity of $f_0$ cannot split such that two or more
$\lambda$-singularities belong to the same fiber.
 \end{enumerate}
\end{theorem}

\begin{remark}\label{r:examplelambda}
Part (a) means that a $\lambda$-singularities cannot completely disappear when deforming
 $f_0$, but of course it may split or its type may change, see \S \ref{examples}.
 The case which is not covered by part (b) can indeed occur, i.e. that some
$\lambda$ splits into $\lambda$'s  along a line
$\{ x(s)\} \times \bC$, see  Example \ref{e1:b}.
\end{remark}

\begin{proof}
  {\bf (a).}  Let $(z,t_0) \in \Sigma_0\times \bC$ be a $\lambda$-singularity
  of $f_0$.  Let us denote by $G(y,s,t)$ the localisation of the map $\tilde
  P(x,x_0,s) -tx_0^d$ at the point $(z,0,t_0) \in \bY$. Let $y_0 =0$ be the
  local equation of the hyperplane at infinity of $\bP^n$. The idea is to
  consider the 2-parameter family of functions $G_{s,t} \colon \bC^{n} \to
  \bC$, where $G_{s,t}(y) = G(y,s,t)$.  Then $G(y,s,t)$ is the germ of a
  deformation of the function $G_{0,t_0}(y)$.

  We consider the germ at $(z,0,t_0)$ of the singular locus $\Gamma$ of the
  map $(G, \sigma, \tau) \colon \bC^{n} \to \bC^3$. This is the union of the
  singular loci of the functions $G_{s,t}$, for varying $s$ and $t$.  We claim
  that $\Gamma$ is a surface, more precisely, that every irreducible component
  $\Gamma_i$ of $\Gamma$ is a surface.  We secondly claim that the projection
  $D\subset \bC^3$ of $\Gamma$ by the map $(y_0,\sigma, \tau)$ is a surface,
  in the sense that all its irreducible components are surfaces. Moreover, the
  projections $\Gamma \stackrel{(y_0,\sigma, \tau)}{\to} D$ and $D
  \stackrel{(s,t)}{\to} \bC^2$ are finite (ramified) coverings.

  All our claims follow from the following fact: the local Milnor number
 conserves in deformations of functions.  The function germ $G_{0,t_0}$ with
 Milnor number, say $\mu_0$, deforms into a function $G_{s,t}$ with finitely
 many isolated singularities, and the total Milnor number is conserved, for
 any couple $(s,t)$ close to $(0,t_0)$.

 Let us now remark that the germ at $(z,0,t_0)$ of $\Sigma \times \bC$ is a
 union of components of $\Gamma$ and projects by $(y_0,s,t)$ to the plane $D_0
 := \{y_0 =0\}$ of $\bC^3$.  However, the inclusion $D_0 \subset D$ cannot be
 an equality, by the above argument on the total ``quantity of singularities''
 and since we have a jump $\lambda >0$ at the point of origin $(z,0,t_0)$. So
 there must exist some other components of $D$.  Every such component being a
 surface in $\bC^3$, has to intersect the plane $D_0\subset \bC^3$ along a
 curve.  Therefore, for every point $(s',t')$ of such a curve, the sum of
 Milnor numbers of the function $G$ on the hypersurface $\{ y_0=0, \sigma =
 s', \tau =t' \}$ (where the sum is taken over the singular points that tend
 to the original point $(z,0,t_0)$ when $s'\to 0$) is therefore strictly
 higher than the one computed for a generic point of the plane $D_0$.
 Therefore our claim (a) will be proved if we prove two things:

\noindent
(i). the singularities of $G$ on the hypersurface $\{ y_0=0, \sigma = s',
 \tau =t' \}$ that tend to the original point $(z,0,t_0)$ when $s'\to 0$
  are included into $G =0$, and

\noindent
  (ii). there exists a component $D_i \subset D$ such that $D_i \cap D_0
 \not= D\cap \{ s=0\}$.

  To show (i), let
  $g_k(y,s)$ denote the degree $k$ part of $P$ after localising it at $p$
   and note that $G(y,s,t) = g_d(y,s) + y_0(g_{d-1}(y,s) +\cdots) -ty_0^d$.
 Then observe that the set:
\begin{equation}\label{eq:indep}
 \Gamma \cap \{y_0 =0\} = \{ \frac{\partial g_d}{\partial y} =0, g_{d-1} =0\}
\end{equation}
does not depend on the variable $t$ and its slice by $\{\sigma =s,\tau =t\}$
 consists of finitely many points. These points may fall into two types:  (I).
 points on  $\{g_d = 0\}$, and therefore on $\{G=0\}$, and  (II). points not
 on $\{g_d = 0\}$. We show that type II points do not actually occur.
This is a consequence of our hypothesis on the constancy of generic singularity
 type at infinity, as follows. By choosing a generic $\hat t$ such that $\hat
 t\not\in \Lambda(s)$ for all $s$, and by using the independence on $t$ of the
 set (\ref{eq:indep}), this condition implies that type II points cannot
 collide
 with type I points along the slice $\{ y_0=0, \sigma = s, \tau = \hat t \}$ as
 $s\to 0$. By absurd, if there were collision, then there would exist a
 singularity in the slice $\{ G= 0, y_0=0, \sigma = 0, \tau = \hat t \}$ with
 Milnor number higher than the generic singularity type at infinity.
It then follows that:
\begin{equation}\label{eq:equal}
\Gamma \cap \{y_0 =0\} = \Gamma \cap \{G= y_0 =0\}
\end{equation}
which proves (i). Now observe that the equality (\ref{eq:equal}) also proves (ii), by a similar
reason:  if there were a component $D_i$ such that $D_i \cap D_0 = D\cap \{ s=0\}$ then there
would exist a singularity in the slice $\{ G= 0, y_0=0, \sigma = 0, \tau = \hat t \}$ with Milnor
number higher than
 the generic singularity type at infinity.
Notice that we have in fact proved more, namely:

\noindent
  (ii'). there is no component $D_i\not= D_0$ such that $D_i \cap D_0 = D\cap
 \{ s=0\}$.

 This ends the proof of (a).

\noindent
{\bf (b).}
  Suppose that there were collision of some singularities out of which
two or more $\lambda$-singularities
 are in the same fibre.
 Then there are at least two different points $z_i \not= z_j$ of $\Sigma_s$
 which collide as $s\to 0$. This situation is excluded by the cgst assumption
(Definition \ref{d:type}).
\end{proof}
%%%%%%%%%%%%%%

%%%%%%%%%%%%%%

%%%%%%%%%%%%%%%%%%%%%
%%%%%%%%%%%%%%%%%%
\section{Local conservation and monodromy in $\mu + \lambda$ constant deformations}\label{s:mulambda}
%%%%%%%%%%%%%%%%%%
%%%%%%%%%%%%%%
In \S\ref{examples} we comment a couple of examples where the inequality of Theorem \ref{t:count}
is strict. Here we show that, when imposing the constancy of $\mu + \lambda$,
 this turns into an equality.
%%%%%%%%%%%%%%%%%%%%%%%%%%%%%%%%%%%%%%
\begin{corollary}\label{c:noloss}
Let $P$ be a constant degree deformation inside the B-class such that
 $\mu(s) + \lambda(s)$ is constant. Then:
\begin{enumerate}
\item As $s\to 0$, there cannot be loss of $\mu$ or of $\lambda$ with
 corresponding atypical values tending to infinity.
 \item
$\lambda$ is upper semi-continuous, i.e. $\lambda(s) \le \lambda(0)$.
\item  there is local conservation of $\mu + \lambda$ at any
 $\lambda$-singularity of $f_0$.
\end{enumerate}
\end{corollary}
\begin{proof}
(a). If there is loss of $\mu$ or of $\lambda$, then
this must necessarily
be compensated by increase of $\lambda$ at some singularity at infinity
 of $f_0$. But Theorem \ref{t:count} shows that
 the local $\mu + \lambda$ cannot increase in the limit.

\noindent
(b). is clear since $\mu(s) + \lambda(s)$ is constant and $\mu(s)$ can only
 decrease when $s\to 0$.

\noindent
(c). Global conservation of $\mu + \lambda$ together with local
 semi-continuity
 (by Theorem  \ref{t:count}) imply local conservation.
\end{proof}
%%%%%%%%%%%%%%
\begin{remark}\label{r:rel}
It is interesting to point out that within the class of B-type polynomials
 there is no
 inclusion relation between the properties ``constant generic singularity
 type'' and ``$\mu(s) + \lambda(s)$ constant'', see Examples \ref{e:2},
\ref{e:1}. We shall see in the following that in the F-class the two conditions
 are equivalent because of the relation (\ref{eq:equality}).
\end{remark}

%%%%%%%
\subsection{Rigidity in deformations with constant
 $\mu + \lambda$}\label{ss:rigid}

 For B-type
 polynomials, we have the formula:
\begin{equation}\label{eq:B}
b_{n-1}(G_s) = \mu(s) +\lambda(s) = (-1)^{n-1} (\chi ^{n,d} - 1) -
\sum_{x\in\Sigma_s}\mu_{x,\gen}(s)
 - (-1)^{n-1} \chi^{\infty}(s),
\end{equation}
where
$\chi ^{n,d} = \chi (V^{n,d}_{gen})  =
  n + 1 - \frac{1}{d} \{ 1 + (-1)^n (d-1)^{n+1} \}$
is the Euler characteristic of the smooth hypersurface $V^{n,d}_{gen}$ of degree $d$
in $\bP^n$ and $\chi^{\infty}(s) := \chi (\{ f_d(x,s) = 0 \})$.
We denote by $\mu_{x,\gen}(s)$ the Milnor number of the singularity
of $\bY_{s,t}$ at the point $(x,t)\in \Sigma_s\times \bC$, for a generic value
of $t$.
The change in  $b_{n-1}(G_s)$
can be described in terms of
change in  $\mu_{x,\gen}(s)$ and $\chi^{\infty}(s)$.
Since the latter is not necessarily semi-continuous (cf Examples
\ref{e:2}--\ref{e:1}), we may
expect interesting exchange of data between the two types of contributions.

\begin{proposition}\label{nolabel}
Let $\Delta \chi^\ity$ denote $(-1)^n(\chi^{\infty}(s)- \chi^{\infty}(0))$.
\begin{enumerate}
\item If $\Delta \chi^\ity < 0$ then the deformation is not cgst.
\item If $\Delta \chi^\ity = 0$ and the deformation has constant $\mu +
  \lambda$ then, for all $x\in \Sigma_s$, $\mu_{x,\gen}(s)$ is constant.
\item If $\Delta \chi^\ity > 0$ then the deformation cannot have constant $\mu +
  \lambda$.
\end{enumerate}
\fin
\end{proposition}

For  F-type polynomials, formula (\ref{eq:B}) takes the following form, see also \cite[(2.1)
 and (2.4)]{ST-gendefo}:

\begin{equation}\label{eq:equality}
 \mu(s)+\lambda(s) = (d-1)^n - \sum_{x\in\Sigma_s}\mu_{x,\gen}(s)-\sum_{x\in
 W_s}\mu_x^\infty(s),
\end{equation}
where $\mu_x^\infty(s)$ denotes the Milnor number of the singularity of $\bY_{s,t}\cap H^\ity$ at
the point $(x,t)\in W_s\times \bC$, which is actually
 independent on the value of $t$. Note that in the F-class we have $\Delta \chi^\ity \ge 0$.

The relation \ref{eq:equality} shows that the change in the Betti number $b_{n-1}(G_s)$ can be
described in terms of change
 in the $\mu_{x,\gen}(s)$
and change in $\mu_x^\infty (s)$. Both are semi-continuous, so they are forced to be constant in
$\mu+\lambda$ constant families.

Consequently, the class of F-type polynomials such that $\mu + \lambda$ is
 constant verifies the hypotheses of Theorem \ref{t:inf}.
It has been noticed by the first named author that in the deformations with
 constant $\mu + \lambda$
 which occur in Siersma-Smeltink's
lists \cite{SS}
 the value of $\lambda$ cannot be dropped to 0.
 Since these deformations are in the F-class and in view of the above observation,
  this behaviour is now completely explained
 by Theorem \ref{t:inf}(a). More precisely, we have proved:
\begin{corollary}\label{c:lambda}
Inside the F-class, a $\lambda$-singularity cannot be
 deformed into only $\mu$-singularities by a
 constant degree deformation with constant $\mu + \lambda$.
\fin
\end{corollary}
%%%%%%%%%%%%%%%%%%%%%%%%%%%%%

%%%%%%%%%%%%%%%%%%%%%%%%%%%%
%%%%%%%%%%%%%%%%%%%%%%%%%%
%%%%%%%%%%%%%%%%%%%%

\subsection{Monodromy in  families with constant $\mu + \lambda$}\label{ss:mon}
%%%%%%%%%%%%%%%%%%%%
For some polynomial $f_0$, one calls {\em monodromy at infinity} the monodromy around a large
enough disc $D$ containing all the atypical values of $f_0$. The locally trivial fibration above
the boundary $\partial \bar D$ of the disc is called {\em monodromy fibration at infinity}.

The global L\^e-Ramanujam problem consists in showing the constancy of the monodromy
 fibration at infinity in a  family with constant
$\mu + \lambda$. Actually one can state the same problem
for any {\em admissible loop} $\gamma$
 in $\bC$, i.e. a simple loop (homeomorphic to a circle) such that it does not
 contain any atypical value of $f_s$, for all $s$ close enough to $0$.

The second named author proved a L{\^e}-Ramanujam type result for a large class of polynomials,
including the B-class (cf \cite{Ti-cras, Ti-reg}),
 with the supplementary condition that there is no loss of $\mu$ at infinity of type \ref{r:ab}(b).
 This hypothesis can now be removed, due to
our Corollary \ref{c:noloss}(a). Moreover, the same result clearly holds over any admissible loop.
Therefore, by revisiting the statement \cite[Theorem 5.2]{Ti-reg}, we get the folowing more
general one:

\begin{theorem}\label{t:mono}
Let $P$ be a constant degree deformation inside the B-class.
 If $\mu + \lambda$ is constant and $n\not=3$ then:
 \begin{enumerate}
\item the monodromy fibrations over any admissible loop are isotopic in the family.
 \item the monodromy fibrations at infinity  are isotopic in the
 family.
 \end{enumerate}
\fin
\end{theorem}

%%%%%%%%%%%%%%%%%%%%%%%%%%%%
%%%%%%%%%%%%%%%%%%%%%%%%%%%%

%%%%%%%%%%%%%%%%%%%%%%%
%%%%%%%%%%%%%%%%%%%%%%%%%%%%%%
\section{Proof of Theorem \ref{t:count}}\label{s:proof}

For the proof, we need to define a certain critical locus. First
endow $\bY$ with the coarsest Whitney stratification $\cW$. Note that (unlike the case of a
single polynomial and its attached space $\bX$ treated in \cite{ST}) {\em we do not require here} that $\bY^\infty$
is a union of strata. Let $\Psi := (\sigma, \tau) : \bY \to
 \bC \times \bC$ be the
projection. The {\em critical locus} $\Crit\hspace{2pt}\Psi$  is the locus
 of points where the restriction of $\Psi$ to some stratum of $\cW$ is not a
submersion.
%%%%%%%%%%%%%%%%%%%%%%%
When writing $\Crit\hspace{2pt}\Psi$ we usually understand a small
 representative
 of the germ of $\Crit\hspace{2pt}\Psi$ at $\bY_{0,*}$.
It follows that $\Crit\hspace{2pt}\Psi$ is a closed analytic set and that its
 affine part $\Crit\hspace{2pt}\Psi  \cap
(\bC^n\times \bC \times \bC)$ is the union, over $s\in \bC$, of the affine critical loci of
the polynomials $f_s$. Notice that both $\Crit\hspace{2pt}\Psi$ and
 its
 affine part $\Crit\hspace{2pt}\Psi
\cap (\bC^n\times \bC\times \bC)$ are in general not product spaces by the
 $t$-variable.
%%%%%%%%%%%%
%%%%%%%%%%%%%%%%%%%%%%%%%%
%\begin{remark}\label{r:B1}
In case of a constant degree one-parametre deformation in the B-class,
the stratification
$\cW$ has a
 maximal stratum which contains the complement of the 2-surface $\Sigma\times \bC$.
 At any point of this complement, all the spaces $\bY$, $\bY_{s,*}$ and
 $\bY_{s,t}$ are nonsingular in the neighbourhood of infinity.
Therefore $\Crit\hspace{2pt}\Psi \cap (\bY^\ity \setminus \Sigma\times \bC)
 = \emptyset$. Since our deformation is in the B-class,
 it follows that the affine part $\Crit\hspace{2pt}\Psi \cap (\bC^n\times \bC\times \bC)$
 is of dimension at most 1. Next, the map is $\Psi$ is submersive over a Zariski-open
subset of any 2-dimensional stratum included in  $\Sigma\times \bC$. It follows
 that the part at infinity of $\Crit\hspace{2pt}\Psi$ has dimension $<2$.
 We altogether conclude that
$\dim\Crit\hspace{2pt}\Psi \le 1$.

%\end{remark}

%%%%%%%%%%%%%%%%%%%%%%%%%%%%%
%\subsection{Proof of Theorem \ref{t:count}}
Nevertheless, this
fact does
not insure that the functions $\sigma$ and $\tau$ have isolated singularity
 with
respect to our stratification
$\cW$. (It is precisely not the case in ``almost all" examples.) Nevertheless,
 in the
pencil
$\sigma + \varepsilon \tau$, $\varepsilon \in \bC$, all the functions except
 finitely
many of them are functions with isolated singularity at $p$ with respect to the stratification
$\cW$. Let us fix some $\varepsilon$ close to zero and consider locally, in some
 good
neighbourhood $\cB$ of $(p,0)\in \bY$,
the couple of functions $\Psi_\varepsilon = (\sigma + \varepsilon \tau, \tau)
\colon \cB \to \bC^2$.

The function $\tau_| \colon (\sigma + \varepsilon \tau)^{-1}(0) \to \bC$
 defines
 an
isolated singularity at $p$ and $(\sigma + \varepsilon \tau)^{-1}(0)$ is a germ of a complete
intersection at $p$. By applying the stratified Bouquet Theorem of \cite{Ti-b}
 we
 get that the Milnor-L{\^e} fibre of $\tau_|$ is homotopy equivalent to a bouquet
 of spheres $\bigvee S^{n-1}$. It follows that the
general fiber of $\Psi_\e$---that is $\cB \cap \Psi_\e^{-1}(s,t)$, for some
 $(s,t)\not\in
\Disc \Psi$---is homotopy equivalent to the same bouquet $\bigvee S^{n-1}$; let $\rho$
 denote
 the number of $S^{n-1}$ spheres in this bouquet.

On the other hand, the Milnor fiber at $p$ of the function $\sigma + \varepsilon
\tau$ is
homotopy equivalent to a bouquet $\bigvee S^{n}$, by the same result {\em
loc.cit.}; let
$\nu$ denote the number of $S^{n}$
spheres.

\begin{figure}[hbtp]
\begin{center}
\epsfxsize=5cm
\leavevmode
\epsffile{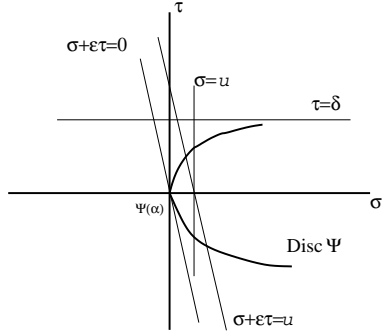}
\end{center}
\caption{{\em
 Counting vanishing cycles.}}
\label{f:1}
\end{figure}

In the remainder, we count the vanishing cycles (I): along $(\sigma +
\varepsilon \tau)^{-1}(0)$, respectively (II): along $(\sigma + \varepsilon
 \tau)^{-1}(u)$, for $u\not= 0$
close enough to $0$, and we compare the results. The vanishing cycles are all in
dimension $n-2$. One may use  Figure \ref{f:1} in order to follow the
 computations; in this picture, the germ
of the  discriminant locus $\Disc \Psi$ at $\Psi(p)$ is the union of the
 $\tau$-axis, $\sigma$-axis  and some other curves.

\noindent
(I). We start with the fiber $\cB \cap \Psi_\e^{-1} (0,\delta)$, where
 $\delta$ is
close
enough to $0$. To obtain $\cB \cap (\sigma + \varepsilon \tau)^{-1}(0)$, which
 is
contractible, one attaches to  $\cB \cap \Psi_\e^{-1} (0,\delta)$ a certain
 number of
$(n-1)$
cells corresponding to the vanishing cycles at infinity, as $t \to 0$, in the
 family
of fibers
$\Psi_\e^{-1} (0,t)$. This is exactly the number $\rho$ defined above and it
 is here the sum of two numbers, corresponding to the attaching in two steps,
as we detail in the following. One is the number of
 cycles in $\cB
\cap \Psi_\e^{-1} (s,\delta)$, vanishing, as $s\to 0$, at points that tend to
$p$ when
$\delta$ tends to $0$; we denote this number by $\xi$. The other number is the
 number of cycles in
$\cB \cap \Psi^{-1}(0,t)$, vanishing as $t\to 0$; this number is
 $\lambda_p(0)$, by definition.
 From this one may draw the inequality: $\lambda_p(0) \le \rho$.

\noindent
(II). Here we start with the fiber $\cB \cap \Psi_\e^{-1} (u,\delta)$, which
 is
homeomorphic
to $\cB \cap \Psi_\e^{-1} (0,\delta)$ and to $\cB \cap \Psi^{-1} (u,\delta)$.
 The
Milnor fiber
$\cB \cap \{\sigma+\e \tau =u\}$  cuts the critical locus $\Crit \Psi$ at
 certain  points
$p_k$. The number of points, counted with multiplicities,
is equal to the local intersection number $\int_p(\{\sigma+\e \tau =0\}, \Crit
\Psi)$.
When walking along $\cB \cap \{\sigma+\e \tau =u\}$, one has to add
to the fiber $\cB \cap \Psi_\e^{-1} (u,\delta)$ a number of cells
corresponding to the
vanishing cycles at points $\{\sigma+\e \tau =u\}\cap \{ \sigma =0\}$,
 which is just the
number $\xi$ defined above, and to the vanishing cycles at points
 $\{\sigma+\e \tau =u\}\cap
\overline{\Crit \Psi \setminus \{ \sigma =0\}}$. The intersection number
$\int_p(\{\sigma+\e \tau =0\}, \overline{\Crit \Psi \setminus \{ \sigma =0\}})$
 is less or equal to the intersection number $\int_p(\{\sigma =0\}, \overline{\Crit
\Psi \setminus \{ \sigma =0\}})$. Now, when walking along $\cB \cap \{\sigma =u\}$, one has to add
to $\cB \cap \Psi^{-1} (u,\delta)$ a number of cells corresponding to the vanishing cycles at
points $p_i$ and $p_j$, which number is, by definition,  $\sum_i \lambda_{p_i}(u) + \sum_j
\mu_{p_j}(u)$. We get the inequality: $\xi + \sum_i \lambda_{p_i}(u) + \sum_j \mu_{p_j}(u) \ge
\rho + \nu$.

 Finally, by collecting the inequalities obtained at steps (I) and (II), we
obtain:
 \begin{equation}\label{eq:lambda}
 \lambda_p(0) = \rho - \xi \le \rho + \nu - \xi  \le \sum_i
\lambda_{p_i}(u) + \sum_j \mu_{p_j}(u),
\end{equation}
 which proves our claim.
 \fin

%%%%%%%%%%%%%%%%%%%%%%%%%%%%%%%%%%%%%%%%%%%%%%
%%%%%%%%%%%%%%%%%%%%%%%%%%%%%%%%%%%%

%%%%%%%%%%%%%%%%%%%%%

%%%%%%%%%%%%%%%%%%%%%%%%%%%%
%%%%%%%%%%%%%%%%%%%%%%%%%%%%
%%%%%%%%%%%%%%%%%%%%%%%%%%%%

\section{Examples}\label{examples}

%%%%%%%%%%%%%%%%%%%%%%%%%%%%

\subsection{F-class examples; behaviour of $\lambda$}

\begin{figure}[hbtp]
\begin{center}
\epsfxsize=14cm
\leavevmode
\epsffile{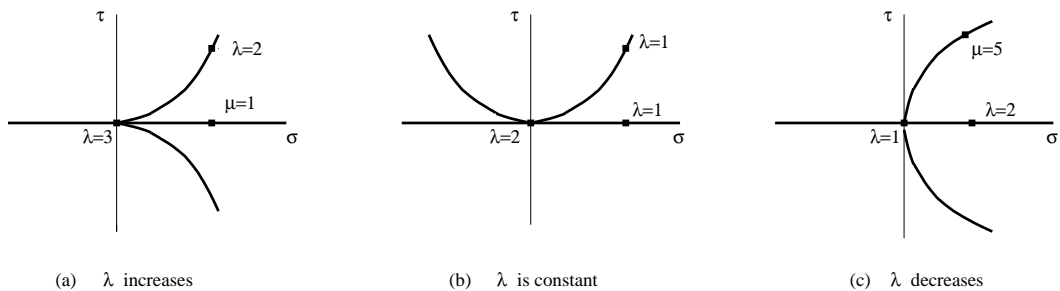}
\end{center}
\caption{{\em Mixed splitting in (a) and (c); pure $\lambda$-splitting in (b).}}
\label{f:2}
\end{figure}

%%%%%%%%%%%%%%%%%%%%%%%%%%%%%%
\begin{example}\label{e1:a}
 $f_s = (xy)^3 + s xy + x$, see Figure \ref{f:2}(a).\\
 This is a deformation inside the F-class, with constant $\mu + \lambda$,
 where $\lambda$ increases.
 For $s \ne 0$: $\lambda = 1+1$  and $\mu = 1$. For $s =  0$:
 $\lambda = 3$ and  $\mu = 0$.
\end{example}
%%%%%%%%%%%%%%%%%%%%%%%%%%%%%%

\begin{example}\label{e1:b}
$f_s=(xy)^4 + s (xy)^2 + x$, see Figure \ref{f:2}(b). \\
  This deformation  has constant $\mu = 0$, $\lambda(0) = 2$ at one point and
 $\lambda(s) = 1+1$  at two points at infinity which differ by the value of
     $t$ only, namely $([0:1], s, 0)$ and $([0:1],s,
  -s^2/4)$.
\end{example}
%%%%%%%%%%%%%%%%%%%%%%%%%%%%%%

\begin{example}\label{e1:c} $f_s = xy^4 + s (xy)^2 + y$, see Figure \ref{f:2}(c). \\
Here  $\lambda$ decreases. For $s \ne 0$: $\lambda = 2$  and $\mu = 5$. For $s =  0$:
 $\lambda = 1$ and  $\mu = 0$.
\end{example}

%%%%%%%%%%%%%%%%%%%%%%%%%%%%%%
%%%%%%%%%%%%%%%%%%%%%%%%%%%%%%

\subsection{B-class examples}

We use in this section formula (\ref{eq:B}). We pay special attention to the
sign of $\Delta \chi^\ity$ and illustrate the difference between cgst-type
deformations and $(\mu+\lambda)$-constant deformations.
%%%%%%%%%%%%%%%%%%%%%%%%%%%%%%

\begin{example}\label{e:2}
$f_{s} = x^4 + sz^4 + z^3 + y$.\\
This is a deformation inside the B-class with constant $\mu + \lambda$,
 which is not cgst at infinity (Definition \ref{d:type}).
We have $\lambda = \mu = 0$ for all $s$. Next, $\bY_{s,t}$ is singular only at $p:=[0:1:0]$ and
the singularities of $\bY_{0,t}^\ity$ change from a single smooth line $\{ x^4=0\}$ with a special
point $p$ on it into the isolated point $p$ which is a $\tilde E_7$ singularity of
$\bY_{s,t}^\ity$. We use the notation  $\oplus$ for the Thom-Sebastiani sum of two types of
singularities in separate variables. We have:

\noindent
$s = 0$:  the generic type is  $A_3 \oplus E_7$ with $\mu= 21$ and
 $\chi(\bY_{0,t}^\ity) = 2$.

\noindent
$s \ne 0$: the generic type is $A_3 \oplus E_6$  with $\mu= 18$ and
$\chi(\bY_{s,t}^\ity) = 5$.

 The jumps of $+3$ and $-3$ compensate each other.
\end{example}

%%%%%%%%%%%%%%%%%%%%%%%%%%%%%%%%%%%%%%%%%%%%%%%%%%%%

\begin{example}\label{e:2a}
$f_{s} = x^4 + sz^4 + z^2 y + z $.\\
This is a $\mu + \lambda$ constant B-type family, with two different singular points of
$\bY_{0,t}$ at infinity, and where the change in one point interacts with the other. It is locally
cgst in one point, but not in the other. We have that $\lambda =3$ and  $\mu = 0$ for all $s$,
$\bY_{s,t}$ is singular  at $p:= [0:1:0]\in H^\ity$ for all $s$ (see types below)
 and at $q:= [1:0:0]\in H^\ity$ with type $A_3$. The singularities of $\bY_{s,t}^\ity$
 change from a single smooth line $\{
x^4=0\}$ into the isolated point $p$ with $\tilde E_7$ singularity.

For the point $p$ we have for all $s$ the
generic type  $A_3 \oplus D_5$  if $t \ne 0$,
which jumps  to $A_3 \oplus D_6$  if $t = 0$.
This causes $\lambda = 3$.

At $q$, the $A_3$-singularity for $s=0$ gets smoothed (independently of $t$) and here the
deformation is not locally cgst. The change on the level of $\chi(\bY_{s,t}^\ity)$ is from 2 to 5,
so $\Delta \chi^\ity =-3$, which compensates the disappearance of the $A_3$-singularity from
 $\bY_{0,t}$ to $\bY_{s,t}$.

% So the contribution in the point $q$ jumped to the point $p$,
% due to the existence
%of the ``channel'' which consists of the line of
%non-isolated singularities $\{ x^4=0\}$.

\end{example}

%%%%%%%%%%%%%%%%%%%%%%%%%%%%%%
\begin{example}\label{e:1}
$f_s = x^2y + x + z^2 + s z^3$.\\
This is a cgst B-type family,
where $\mu + \lambda$ is not constant. Notice that
$f_s$ is F-type for all $s\not= 0$, whereas $f_0$ is not F-type (but still B-type).
The generic type at infinity is $D_4$ for all $s$ and there is a jump
 $D_4 \to D_5$ for $t= 0$ and all $s$.
For $s \ne 0$ a second jump $D_4 \to D_5$ occurs for  $t= c/s^2$, for some
constant $c$.

There are no affine critical points,
i.e. $\mu(s) = 0$ for all $s$, but
$\lambda(s) = 2$ if $s \ne 0$ and $\lambda(0) =1$.
We have that $\Lambda(f_s) = \{ 0, c/s^2\}$ for all $s\not=0$,
and that $\chi^{\infty}$ changes from $3$ if $s=0$ to $2$ if $s\not= 0$, so
$\Delta\chi^\ity = +1$.

There is a persistent $\lambda$-singularity in the fibre over $t=0$
and there is a branch of the critical locus $\Crit \Psi$ which is asymptotic
to $t=\infty$.
\end{example}
%%%%%%%%%%%%%%%%%%%%%%%%%%%%%%
%%%%%%%%%%%%%%%%%%%%%%%%%%%%%%
\subsection{Cases of lower semi-continuity at $\lambda$-singularities}\label{ss:ineq}

In Theorem \ref{t:count} we have an inequality which we may write in short-hand as follows, by
referring to its proof (formula \ref{eq:lambda}):
\begin{equation}\label{eq:lambda2}
\lambda = I_{gen} - \nu \le I_{gen} \le I_{s=0}
\end{equation}
This inequality can have two different sources:
\begin{itemize}
\item[-] the nongeneric intersection number $I_{s=0}$ and its difference to
  the generic one $I_{gen}$,
\item[-] the number $\nu$, which is related to the equisingularity
properties of $\bY$.
\end{itemize}

So the excess in the formula is $\nu + (I_{s=0} - I_{gen})$.
 The following examples illustrate the different types of excess:
$\nu\not= 0$,
respectively $\nu= 0$ and $I_{s=0} - I_{gen}>0$. In the latter case, the space
 $\bY$ is singular.

\begin{example}\label{e:excess1}
We start with a F-type polynomial $f_0$ and consider a Yomdin deformation $f_0 - s x_1^d$ for
sufficient general $x_1$.  In this case the space $\bY$ is non-singular and the function $\sigma +
\varepsilon \tau$ behaves locally as a linear function.  It follows that $\nu = 0 $. Moreover in
this case $I_{s=0} - I_{gen}$ turns out to be positive because of the tangency of some components
of the discriminant set to the $s$-axis.  Compare to \cite[Theorem 5.4]{ST-gendefo}, where the
local lower semi-continuity was proved in the case of Yomdin deformations.

\end{example}

\begin{example}\label{e:excess2}
$f_s = x ^2y ^b + x + s x y ^k $.\\ In the range $\frac{b}{2} < k \le b$, this has the following
data:\\ $s=0$: $\lambda = b$, $\mu = 0$, $\lambda + \mu = b$;\\ $s \ne 0$: $\lambda = 0$, $\mu =
2k$, $\lambda + \mu = 2k$.\\ Both intersection numbers $I_{gen}$ and $I_{s=0}$ are the same and
equal to $2k$. We read the inequality (\ref{eq:lambda2}) as: $b = 2k - \nu \le 2k \le 2k$.  So
$\nu = 2k - b $ and this is possitive in case $\frac{b}{2} < k \le b$.

For the complementary range $1 < k < \frac{b}{2}$ we have a family with an
extra $\lambda$-discriminant branch at $t=0$.
There is the following data here:\\
$s=0$: $\lambda = b$, $\mu = 0$, $\lambda + \mu = b$ ;\\ $s \ne 0$: $\lambda =
b - 2k$, $\mu = 2k$, $\lambda + \mu = b$.\\ In this range one has $\nu = 0$,
$\lambda = b = I_{s=0}$, which gives equality in Theorem \ref{t:count}. This
local conservation is characteristic to families with constant global
$\mu+\lambda$, see Corollary \ref{c:noloss}(c).
\end{example}

%%%%%%%%%%%%%%%%%%%%%%%%%%%%

%%%%%%%%%%%%%%%%%%%%%%%%%%

%%%%%%%%%%%%%%%%%%%%%%%%%%%
%%%%%%%%%%%%%%%%%%%%%%%%%%%
%%%%%%%%%%%%%%%%%%%%%%%%%%%

%%%%%%%%%%%%%%%%%%
%%%%%%%%%%%%%%%%%%
%%%%%%%%%%%%%%%%%%
%%%%%%%%%%%%%%%
%

%%%%%%%%%%%%%%%%%%%%%%%%%%%%
%%%%%%%%%%%%%%%%%%%%%%%%%%%
%%%%%%%%%%%%%%%%%%%%%%%%%%%%%%%%%%%%%%%%%%%%%%%%%%%%%%

%%%%%%%%%%%%%%%%%%


\begin{thebibliography}{MMM}
\footnotesize{

\bibitem[AC]{AC}
 N. A'Campo,  {\em Le nombre de Lefschetz d'une monodromie},
 Indag. Math.  35 (1973),  113--118.

\bibitem[AF]{AF}
A. Andreotti, T. Frankel, {\em The Lefschetz theorem on hyperplane sections},
 Ann. of Math.
(2)  69 (1959), 713--717.


 \bibitem[Br]{Br}
  S.A. Broughton, {\em On the topology of polynomial
hypersurfaces},
Proceedings A.M.S. Symp. in Pure. Math., vol.  40, I (1983),
165-178.

\bibitem[La]{La}
F.  Lazzeri, 
{\em  A theorem on the monodromy of isolated singularities}, in:
 Singularit\'es \`a Carg\`ese, 1972,  pp. 269--275. Asterisque, Nos. 7 et 8,
 Soc. Math. France, Paris, 1973.

 \bibitem[L\^e]{Le}
 L\^e D.T.,  {\em  Une application d'un
 th\'eor\`eme d'A'Campo
\`a l'\'equisingularit\'e}, Nederl. Akad. Wetensch. Proc. (Indag. Math.)
35 (1973),  403--409.

\bibitem[NS]{NS}
A. N{\'e}methi, C. Sabbah,  {\em Semicontinuity of the spectrum at infinity}, Abh. Math. Sem. Univ.
Hamburg 69 (1999), 25--35.

\bibitem[Pa]{Pa}
 A. Parusi\'nski, {\em On the bifurcation set of complex polynomial with
 isolated
singularities at infinity}, Compositio Math. 97 (1995), no. 3, 369--384.


\bibitem[SS]{SS}
D. Siersma, J. Smeltink,
       {\em  Classification of singularities at infinity of polynomials of
degree $4$ in two variabales},  Georgian Math. J.  7 (1)
(2000), 179--190.

\bibitem[ST1]{ST}
 D. Siersma, M. Tib\u ar, {\em Singularities at infinity
and their
vanishing cycles},  Duke Math. Journal 80 (3) (1995), 771-783.


\bibitem[ST2]{ST-gendefo}
D. Siersma, M. Tib\u ar, {\em Deformations of polynomials, boundary
  singularities
 and monodromy}, Mosc. Math. J. 3 (2) (2003), 661--679.

\bibitem[Ti1]{Ti-b}
M. Tib\u ar, {\em Bouquet decomposition of the Milnor fibre},
 Topology  35, 1 (1996), 227--241.

\bibitem[Ti2]{Ti-cras}
M. Tib\u ar, {\em On the monodromy fibration of polynomial functions with
singularities at infinity},
 C.R. Acad. Sci. Paris,  324 (1997), 1031--1035.

\bibitem[Ti3]{Ti-reg}
M. Tib\u ar, {\em Regularity at infinity of real and complex polynomial maps},
 in: Singularity
Theory, The C.T.C Wall Anniversary Volume,
LMS Lecture Notes Series 263 (1999), 249--264. Cambridge
University Press.


 }
\end{thebibliography}
\end{document}